\documentclass[12pt,letterpaper,reqno]{amsart}

\usepackage{times}
\usepackage[T1]{fontenc}
\usepackage{mathrsfs}
\usepackage{latexsym}
\usepackage[dvips]{graphics}
\usepackage{epsfig}
\usepackage{hyperref}

\usepackage{amsmath,amsfonts,amsthm,amssymb,amscd}
\input amssym.def
\input amssym.tex

\newcommand{\kommentar}[1]{}

\newcommand\be{\begin{equation}}
\newcommand\ee{\end{equation}}
\newcommand\bea{\begin{eqnarray}}
\newcommand\eea{\end{eqnarray}}
\newcommand\bi{\begin{itemize}}
\newcommand\ei{\end{itemize}}
\newcommand\ben{\begin{enumerate}}
\newcommand\een{\end{enumerate}}
\newcommand\bc{\begin{center}}
\newcommand\ec{\end{center}}
\newcommand\ba{\begin{array}}
\newcommand\ea{\end{array}}

\theoremstyle{definition}

\newcommand{\erf}{{\rm erf}}

\numberwithin{equation}{section}

\begin{document}

\title{The `Real' Schwarz Lemma}

\author{Steven J. Miller}\email{Steven.J.Miller@williams.edu}\address{Department of Mathematics and Statistics, Williams College, Williamstown, MA 01267}

\author{David A. Thompson}\email{David.A.Thompson@williams.edu}\address{Department of Mathematics and Statistics, Williams College, Williamstown, MA 01267}




\date{\today}

\thanks{We thank our classmates from Math 302: Complex Analysis (Williams College, Fall 2010) for many enlightening conversations, especially David Gold and Liyang Zhang. The first named author was partially supported by NSF grant DMS0970067. }

\begin{abstract} The purpose of this note is to discuss the real analogue of the Schwarz lemma from complex analysis. \end{abstract}


\maketitle

\setcounter{equation}{0}



One of the most common themes in any complex analysis course is how different functions of a complex variable are from functions of a real variable. The differences can be striking, ranging from the fact that any function which is complex differentiable once must be complex differentiable infinitely often \emph{and} further must equal its Taylor series, to the fact that any complex differentiable function which is bounded must be constant. Both statements fail in the real case; for the first consider $x^3 \sin(1/x)$ while for the second just consider $\sin x$. In this note we explore the differences between the real and complex cases of the Schwarz lemma.

The Schwarz lemma states that if $f$ is a holomorphic map of the unit disk to itself that fixes the origin, then $|f'(0)| \le 1$; further, if $|f'(0)| = 1$ then $f$ is an automorphism (in fact, a rotation). What this means is that we cannot have $f$ locally expanding near the origin in the unit disk faster than the identity function, even if we were willing to pay for this by having $f$ contracting a bit near the boundary. The largest possible value for the derivative at the origin of such an automorphism is 1. This result can be found in every good complex analysis book (see for example \cite{Al,La,SS}), and serves as one of the key ingredients in the proof of the Riemann Mapping Theorem. For more information about the lemma and its applications, see the recent article in the Monthly by Harold Boas \cite{Bo}.


It's interesting to consider the real analogue. In that situation, we're seeking a real analytic map $g$ from $(-1,1)$ to itself with $g(0) = 0$ and derivative $g'(0)$ as large as possible. After a little exploration, we
quickly find two functions with derivative greater than 1 at the origin. The first is $g(x) =\sin(\pi x / 2)$, which has $g'(0) = \pi/2 \in (1,2)$. The second is actually an infinite family: letting $g_a(x) = (a+1)x/(1+ax^2)$ we see that $g_a$ is real analytic on $(-1,1)$ so long as $|a| \le 1$, and $g_a'(0) = 1+a$. Using this example, we see we can get the derivative as large as 2 at the origin. Unfortunately, if we take $|a| > 1$ then $g_a$ is no longer a map from $(-1,1)$ to itself; for example, $g_{1.01}(.995) > 1.00001$.

Notice both examples fail if we try to extend these automorphisms to maps on the unit disk. For example, when $z = 3 i / 5$ then already $\sin(\pi z/2)$ has absolute value exceeding 1, and thus we would not have an automorphism of the disk. For the family $g_a$, without loss of generality take $a > 0$. As $z \to i$ then $g_a(z) \to \frac{1+a}{1-a}i$, which is outside the unit disk if $a > 0$.


While it is easy to generalize our family $\{g_a\}$ to get a larger derivative at 0, unfortunately all the examples we tried were no longer real analytic on the entire interval $(-1,1)$. As every holomorphic function is also analytic (which means it equals its Taylor series expansion), it seems only fair to require this property to hold in the real case as well. Interestingly, there is a family of real analytic automorphisms of the unit interval fixing the origin whose derivatives become
arbitrarily large at 0. Consider $h_k(x) = \erf(kx) / \erf(k)$, where $\erf$ is the error function: \be \erf(x) \ := \ \frac{2}{\sqrt{\pi}} \int_0^x e^{-t^2} dt. \nonumber\ \ee The error function has a series expansion converging for all complex $z$, \be \erf(z) \ = \ \frac2{\sqrt{\pi}} \sum_{n=0}^\infty \frac{(-1)^n z^{2n+1}}{n!(2n+1)} \ = \ \frac2{\sqrt{\pi}} \left(z - \frac{z^3}{3} + \frac{z^5}{10} - \frac{z^7}{42} + \cdots \right)\nonumber\ \ee (this follows by using the series expansion for the exponential function and interchanging the sum and the integral),  and is simply twice the area under a normal distribution with mean 0 and variance 1/2 from $0$ to $x$. From its definition, we see $\erf(-x) = -\erf(x)$, the error function is one-to-one, and for $x \in (-1,1)$ our function $\erf(kx)/\erf(k)$ is onto $(-1,1)$.

\begin{figure}
\begin{center}
\scalebox{.6}{\includegraphics{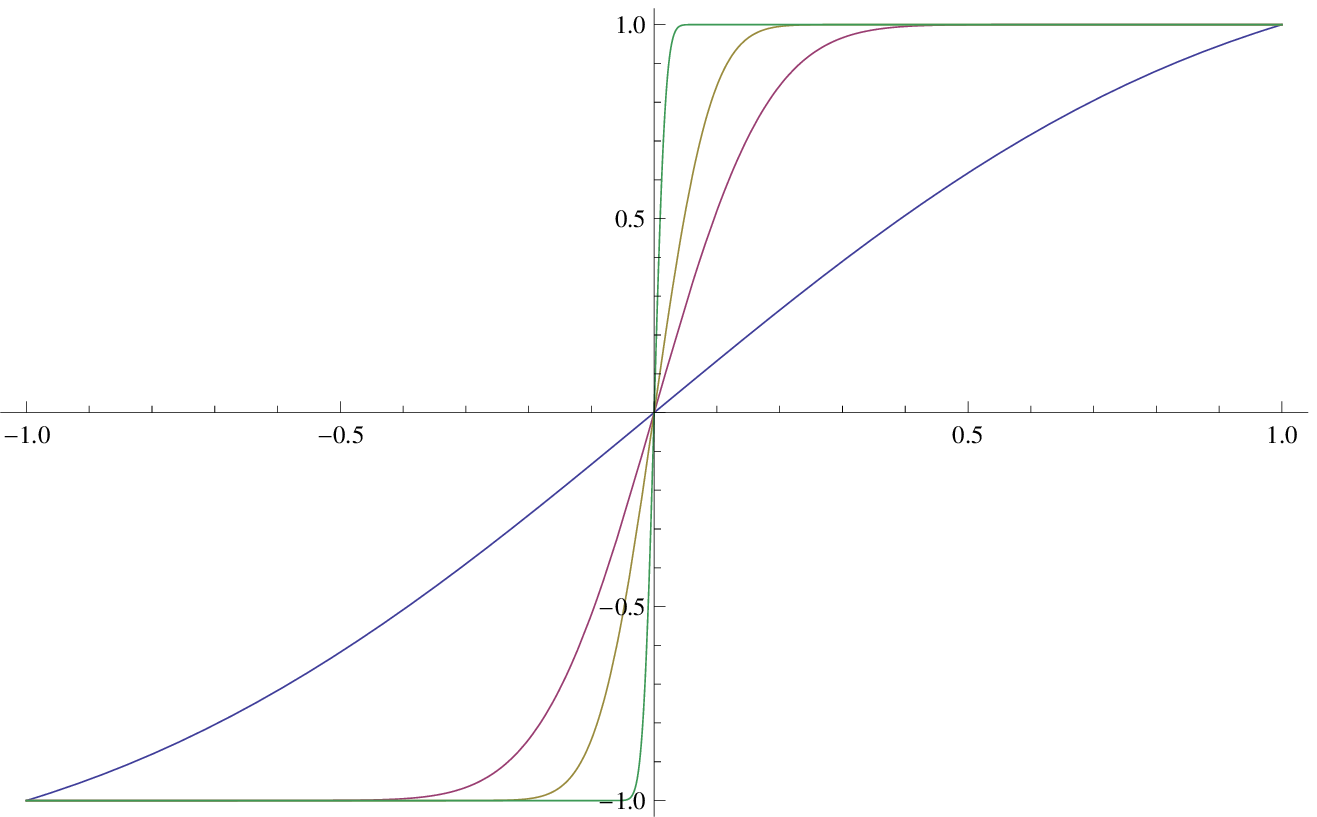}}\ \scalebox{.6}{\includegraphics{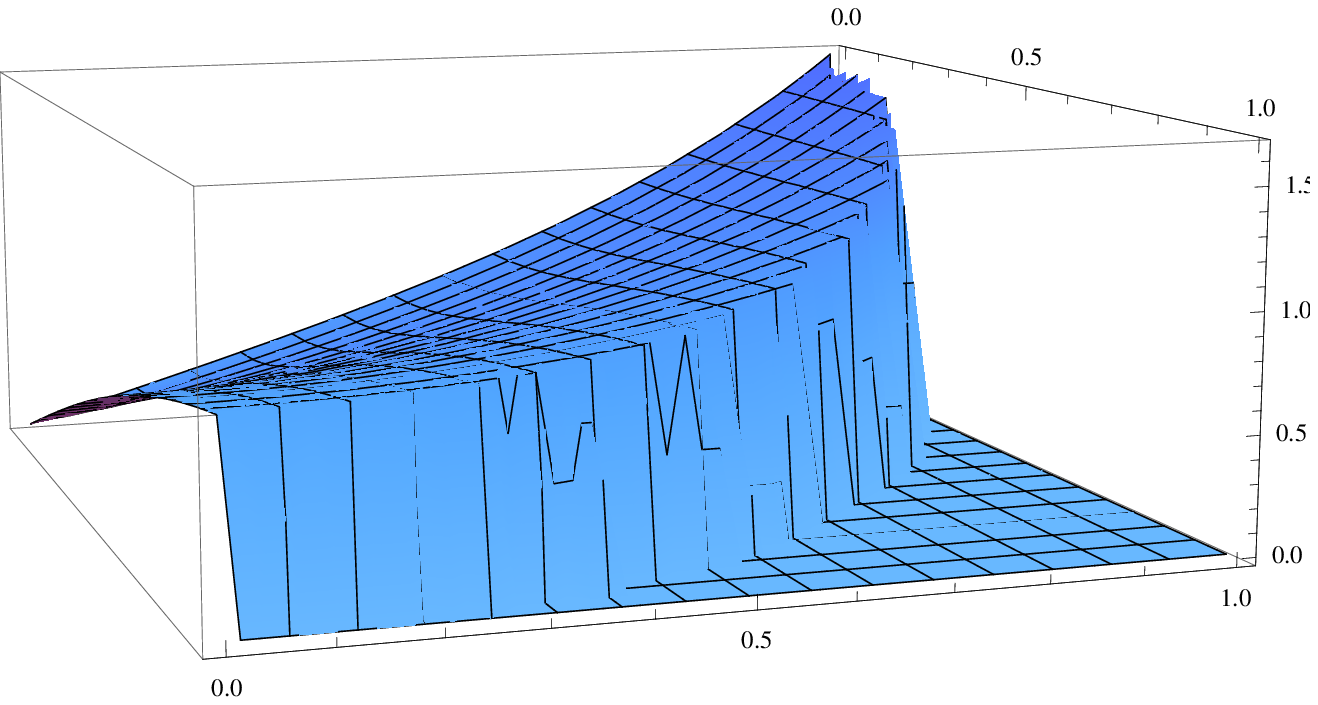}}\
\caption{\label{fig:erf} Plot of the scaled error functions. (1) Left: ${\rm Erf}(kx)/{\rm Erf}(x)$ for $k \in \{1, 5, 10, 50\}$ and $x \in (-1,1)$; (2) Right: Plot of $|{\rm Erf}(z)|$ for $|z| \le 1$.}
\end{center}\end{figure}

Using the Fundamental Theorem of Calculus, we see that $h_k'(x) = 2 \exp(-k^2 x^2) k) / \sqrt{\pi} \erf(k)$, and thus $h_k'(0) = 2k/\sqrt{\pi} \erf(k)$. As $\erf(k) \to 1$ as $k\to\infty$, we find $h_k'(0) \sim 2k/\sqrt{\pi} \to \infty$, which shows that, yet again, the real case behaves in a markedly different manner than the complex one. As $h_k$ is an entire function with large derivative at 0, if we regard it as a map from the unit disk it must violate one of the conditions of the Schwarz lemma. From the series expansion of the error function, it's clear that if we take $z = iy$ then $h_k(iy)$ tends to infinity as $y\to 1$ and $k\to\infty$; thus $h_k$ does not map the unit disk into itself, and cannot be a conformal automorphism (see Figure \ref{fig:erf} for plots in the real and complex cases).

\ \\


\begin{thebibliography}{BCD}


\bibitem[Al]{Al}
L. Ahlfors, \emph{Complex Analysis}, 3rd edition, McGraw-Hill, New
York, 1979.

\bibitem[Bo]{Bo}
H. P. Boas, \emph{Julius and Julia: Mastering the Art of the Schwarz Lemma}, American Mathematical Monthly \textbf{117} (2010) 770--785.


\bibitem[La]{La}
S. Lang, \emph{Complex Analysis}, Graduate Texts in Mathematics,
Vol. 103, Springer-Verlag, New York, 1999.


\bibitem[SS]{SS}
E. Stein and R. Shakarchi, \emph{Complex Analysis}, Princeton
University Press, Princeton, NJ, 2003.




\end{thebibliography}
\end{document}